\documentclass[12pt]{amsart}
\usepackage{amsfonts}
\usepackage{amssymb}
\theoremstyle{definition}

\theoremstyle{remark}

\newcommand{\R}{\mathbb R}

\newcommand{\ds}{\displaystyle}

\begin{document}

\centerline{\bf ANNUAIRE \ DE \ L'UNIVERSIT\'E \ DE \ SOFIA \ "KLIMENT \ OHRIDSKI"}

\vspace{0.4in}
\centerline{\bf FACULT\'E DE MATH\'EMATIQUES ET M\'ECANIQUE}

\vspace{0.4in}
\centerline{\it Tome 81  \qquad\qquad\qquad\qquad\qquad Livre 1 -- Math\'ematiques   \qquad\qquad\qquad\qquad\qquad   1987}

\vspace{0.8in}

\centerline{\large\bf CONGRUENCE OF HYPERSURFACES  }
\centerline{\large\bf OF A PSEUDO-EUCLIDEAN SPACE}

\vspace{0.6in}
\centerline{\bf  Ognian Kassabov}

\vspace{0.4in}
{\sl  R. S. Kulkarni has proved [1] that the so-called bending of a hypersurface in an
Euclidean space determines the congruence class of the hypersurface. In the present 
paper we show that a similar result holds for hypersurfaces of a pseudo-Euclidean 
space $\R_s^{n+1}$, $n>2$. We prove also a corresponding theorem, which accounts for 
the behaviour of the second fundamental form on isotropic vectors.}

\vspace{0.5in}
\centerline{\bf 1. PRELIMINARIES }

\vspace{0.3in}
Let $M$ be a Riemannian or a pseudo-Riemannian manifold with a metric tensor $g$.
A tangent vector $\xi$ is said to be isotropic, it it is nonzero and $g(\xi,\xi)=0$.
Of course, for isotropic vectors one speaks only when the manifold is pseudo-Riemannian,
i.e. when $g$ is an indefinite metric. The values of a symmetric tensor of type (0,2)
on isotropic vectors give a good information about this tensor, as it is shown by
the following

{\bf Lemma 1} [2]. Let $M$ be a pseudo-Riemannian manifold. If $L$ is a symmetric
tensor of type (0,2) on a tangent space $T_pM$, such that $L(\xi,\xi)=0$ for 
every isotropic vector $\xi$ in $T_pM$, then $L=cg$, where $c$ is a real number.

Let $\nabla$ and $R$ denote the covariant differentiation and the curvature tensor ot
$M$, respectively. The Ricci tensor and the scalar curvature will be denoted by $S$
and $\tau$, respectively. Then the Weil conformal curvature tensor $C$ for $M$ is
defined by
$$
	C=R-\frac{1}{n-2}\varphi+\frac{\tau}{(n-1)(n-2)}\pi_1 \ ,
$$
where $n=$\,dim\,$M$, $\varphi$ is defined by
$$
	\varphi (T)(x,y,z,u)=g(x,u)T(y,z)-g(x,z)T(y,u)+g(y,z)T(x,u)-g(y,u)T(x,z)\ ,
$$
for any symmetric tensor $T$ of type (0,2) and $\pi_1=\frac12\varphi(g)$. As it
is well known [3], if $n>3$, then $M$ is conformally flat if and only is the
Weil conformal curvature tensor vanishes identivally. If $n=3$ a necessary and
sufficient condition for $M$ to be conformally flat is [3]
$$
	\left(\nabla_X\left(S-\frac\tau4\,g\right)\right)(Y,Z)-
	  \left(\nabla_Y\left(S-\frac\tau4\,g\right)\right)(X,Z)=0 \ .   \leqno (1.1)
$$

If\, $\overline M$\, is another Riemannian or pseudo-Riemannian manifold, we denote
the corres\-ponding objects for $\overline M$ by a bar overhead. Assume that $f$ 
is a conformal diffeomorphism of $M$ onto $\overline M$: $f^*\bar g=\varepsilon
e^{2\sigma}g$, where $\varepsilon = \pm 1$ and $\sigma$ is a smooth function. 
Then we have [3]
$$
	f^*\overline R=\varepsilon e^{2\sigma} \{R+\varphi(Q) \} \ ,  \leqno(1.2)
$$
where
$$
	Q(X,Y)= X\sigma Y\sigma-g(\nabla_X\nabla\sigma,Y)
				-\frac12\left\Arrowvert \nabla\sigma\right\Arrowvert^2 g(X,Y) \ ,   \leqno(1.3)
$$
$\nabla\sigma$ denoting the gradient of $\sigma$ and 
$\left\Arrowvert \nabla\sigma\right\Arrowvert^2=g(\nabla\sigma,\nabla\sigma)$.

In [4] we have proved the following

{\bf Lemma 2.} Let $M$ and $\overline M$ be pseudo-Riemannian manifolds of dimension 
$>2$ and $f$ be a diffeomorphism of $M$ onto $\overline M$. Assume that at a point $p$ 
of $M$ there exists an isotropic vector $\xi$, such that every isotropic vector, which
is sufficiently close to $\xi$, is mapped by $f_*$ in an isotropic vector in $f(p)$. Then
$f_*$ is a homothety at $p$. 

In what follows $M$ will be a hypersurface of an Euclidean space $\R^{n+1}$ or of a
pseudo-Euclidean space $\R_s^{n+1}$, such that the restriction $g$ of the usual metric
of $\R_s^{n+1}$ to $M$ is nondegenerate. Denote the second fundamental form of $M$ by
$h$. Then we have the equation of Gauss
$$
	R(X,Y,Z,U)=h(X,U)h(Y,Z)-h(X,Z)h(Y,U)
$$
and the equation of Codazzi
$$
	(\nabla_Xh)(Y,Z)-(\nabla_Yh)(X,Z)=0 \ .
$$
Recall also, that a point $p$ of $M$ is said to be quasi-umbilic, if
$$
	h=\alpha\, g + \beta\,\omega\otimes \omega
$$
in $p$, where $\alpha,\,\beta$ are real functions and $\omega$ is an 1-form.
In particular, if $\beta$ is zero the point $p$ is called umbilic.

The bending [1] $K_h$ of $M$ is said to be  the function, assigning to each nonisotropic
nonzero tangent vector $x$ at a point of $M$ the number
$$
	K_h(x)=\frac{h(x,x)}{g(x,x)} \ .
$$ 
Two hypersurfaces $M$ and $\overline M$ being defined a diffeomorphism $f$ of $M$
onto $\overline M$ is said to be bending preserving [1], if
$$
	\overline K_{\bar h}(f_*x)=K_h(x)    \leqno (1.4)
$$
for each nonisotropic nonzero vector $x$ on $M$, whose image is also nonisotropic.
The analogue of (1.4) for isotropic vectors is
$$
	\lim_{x\rightarrow\xi}\frac{\overline K_{\bar h}(f_*x)}{K_h(x)} =1 \ ,    \leqno (1.5)
$$
where the isotropic vector $\xi$ is approximated by nonisotropic nonzero vectors, whose
images are also nonisotropic. We shall prove:

{\bf Theorem 1.} Let $M$ and $\overline M$ be hypersurfaces with indefinite metrics
in $\R_s^{n+1}$, $n>2$, and let $f$ be a diffeomorphism of $M$ onto $\overline M$,
satisfying (1.5) for each isotropic vector $\xi$ on $M$. If the nonquasi-umbilic
points are dense in $M$ and if $M$ is not conformally flat, then $f$ is a congruence.

We recall that $f$ is said to be a congruence if it can be extended to a motion
of $\R_s^{n+1}$.

{\bf Theorem 2.} Let $M$ and $\overline M$ be hypersurfaces with indefinite metrics
in $\R_s^{n+1}$, $n>2$, and let $f$ be a bending preserving diffeomorphism of $M$ 
onto $\overline M$. If the nonumbilic points are dense in $M$ and the curvature tensor
of $M$ does not vanish identically in a point $p$, then there exists a neighbourhood
$V$ of $p$ such that $f_{|V} $ is a congruence of $V$ onto $f(V)$.

{\bf Remark.} The proof of the congruence theorem in [1] can be applied for
hypersurfaces with definite metrics (i.e. spacelike hypersurfaces) in $\R_1^{n+1}$.

\vspace{0.5in}

\centerline{\bf 2. BASIC RESULTS}

\vspace{0.3in}
In this section we prove two lemmas, which will be useful in the proofs of
Theorems 1 and 2.

{\bf Lemma 3.} Let $M$ and $\overline M$ be hypersurfaces with indefinite metrics
in $\R_s^{n+1}$, $n>2$, and let $f$ be a diffeomorphism of $M$ onto $\overline M$,
satisfying (1.5) for each isotropic vector $\xi$ on $M$. If the nonumbilic points
are dense in $M$, then

a) $f$ is conformal: $f^*\bar g=\varepsilon e^{2\sigma}g$;

b) $f^*\bar h=\varepsilon e^{2\sigma}\{ h+\lambda g \}$, where $\lambda$ is a smooth function;

c) $f^*\overline R=e^{4\sigma}\{ R+\lambda\varphi(h)+\lambda^2\pi_1 \}$;

d) the following equations hold:
$$
	X\sigma B(Y,Z)-Y\sigma B(X,Z) +\frac1n\{ X\lambda\, g(Y,Z) -Y\lambda\,g(X,Z) \}=0 \ ,
			\leqno (2.1)
$$
$$
	B(Y,\nabla\sigma) = \frac{n-1}{n} \,Y\lambda \ ,       	\leqno (2.2)
$$
where $\ds B=h-\frac{{\rm tr}h}n g$.

{\it Proof.} Let $p$ be a nonumbilic point of $M$, i.e. $h$ is not proportional to
$g$ in $p$. Then by Lemma 1 there exists an isotropic vector $\xi$ in $T_pM$, such that 
$h(\xi,\xi) \ne 0$. Hence $h(\xi',\xi')\ne 0$ for each isotropic vector $\xi$, which
is sufficiently close to $\xi$. Then (1.5) implies that $f_{*p}\xi'$ is isotrpic.
According to Lemma 2, $f_*$ is a homothety at $p$. Since the nonumbilic points are
dense in $M$, $f$ is conformal and then a) is proved.

From a) and (1.5) it follows $(f^*\bar h)(\xi,\xi)=\varepsilon e^{2\sigma}h(\xi,\xi)$.
Applying again Lemma 1, we obtain b). Then c) follows from b) and from the equations of 
Gauss for $M$ and $\overline M$.

To simplify the notations in the proof of d), we identify $M$ with $\overline M$ via $f$,
and omit $f^*$ from the formulas. Then we have [3]
$$
	\overline\nabla_XY=\nabla_XY +X\sigma Y+Y\sigma X - g(X,Y)\nabla\sigma \ .
$$
Hence, using b) and the equations of Codazzi for $g$ and $\bar g$, we find
$$
	\begin{array}{c}
		X\sigma\, h(Y,Z) - Y\sigma\,h(X,Z)+X\lambda\,g(Y,Z)-Y\lambda\, g(X,Z) \\
		+g(X,Z)h(Y,\nabla\sigma)-g(Y,Z)h(X,\nabla\sigma)=0 \ ,
	\end{array}   \leqno (2.3)
$$
which implies immediately
$$
	h(Y,\nabla\sigma)=\frac{n-1}{n}\,Y\lambda + \frac{{\rm tr} \,h}{n}\,Y\sigma \ ,  \leqno(2.4)
$$
i.e. (2.2). From (2.3) and (2.4) we obtain (2.1). This proves the lemma.

We note that the conditions of Lemma 3 are fulfilled in Theorem 1, as well as in Theo\-rem 2.

{\bf Lemma 4.} If in Lemma 3 $\left\Arrowvert\sigma\right\Arrowvert^2=0$ and $U$ denotes the open
set $\{ p\in M\ : \ \nabla\sigma \ne 0 \}$, then

a) each point of $U$ is quasi-umbilic;

b) $R=0$ in $U$.

{\it Proof.} We shall use a connected component $U_1$ of $U$. Let in (2.1) 
$X=Z=\nabla\sigma$, $Y=\nabla\lambda$. By (2.2), we obtain 
$$
	(\nabla\sigma)\lambda = 0 \ .   \leqno (2.5)
$$
Now we put $X=\nabla\sigma$ in (2.1) and get use of (2.2) and (2.5). The result is
$Y\sigma\,Z\lambda=0$ for arbitrary vector fields $Y$, $Z$. Since $\nabla\sigma$ can 
not vanish in $U_1$, it follows $\lambda=$\,const (in $U_1$). Then (2.1) reduces to
$$
	X\sigma\, B(Y,Z)-Y\sigma\,B(X,Z)=0 \ ,
$$
which implies
$$
	B=\mu\,d\sigma \otimes d\sigma \ ,   \leqno (2.6)
$$
where $\mu $ is a smooth function. Equivalently, we may write
$$
	h=\frac{{\rm tr}\,h}{n} g+\mu\,d\sigma \otimes d\sigma \ ,   \leqno (2.6')
$$
thus proving a). From the equation of Gauss for $g$ it follows
$$
	S(x,x)={\rm tr}\,h.h(x,y)-\sum_{i=1}^n h(x,e_i)h(y,e_i)g(e_i,e_i) \ ,
$$
where $\{ e_i;\, i,\hdots,n \}$ is an orthogonal frame. Hence, using $(2.6')$, we obtain
$$
	S=\frac{n-2}n \mu\,{\rm tr}\,h\,d\sigma\otimes d\sigma +\frac{n-2}{n^2}({\rm tr}\,h)^2g \ .
	          \leqno (2.7)
$$
Thus we get
$$
	\tau = \frac{n-1}n({\rm tr}\,h)^2 \ .    \leqno (2.8)
$$
From (2.7) and (2.8) we compute for $\ds P=S-\frac{\tau}n\,g$:
$$
	P = \frac{n-2}n\,\mu\,{\rm tr}\,h\, d\sigma\otimes d\sigma \ . \leqno (2.9)
$$
By Lemma 3 c) we find immediately
$$
	\begin{array}{l}
		f^*\overline S=\varepsilon e^{2\sigma}\{ S+\lambda(n-2)h+\lambda\,{\rm tr}\,h.g
										+(n-1)\lambda^2g\}\ ,  \\
		f^*\bar\tau=\tau+2(n-1)\lambda\,{\rm tr}\,h + n(n-1)\lambda^2\ ,  \\
		f^*\overline P=\varepsilon e^{2\sigma}\{ P+(n-2)\,\lambda\,B \} \ .	
	\end{array}   \leqno (2.10)
$$
Analogously, (1.2) yields
$$
	f^*\overline P = P+(n-2)Q-\frac{n-2}{2n(n-1)}(\varepsilon\bar\tau\,e^{2\sigma}-\tau)g\ .
$$
From the last two equations we obtain
$$
	Q=\frac{\varepsilon e^{2\sigma}-1}{n-2}P+\varepsilon\lambda\, e^{2\sigma}B+
								\frac{\varepsilon\bar\tau e^{2\sigma}-\tau}{2n(n-1)}\,g \ .
$$
Hence, using (2.6) and (2.9), we find
$$
	Q=\nu\,d\sigma\otimes d\sigma + \frac{\varepsilon\bar\tau e^{2\sigma}-\tau}{2n(n-1)}g \ ,  \leqno(2.11)
$$
where
$$
	\nu =\mu\left( \frac{\varepsilon e^{2\sigma}-1}{n}\,{\rm tr}\,h+\varepsilon\lambda e^{2\sigma} \right) \ .
$$
Since $\nabla\sigma$ is isotropic, (1.3) yields $Q(X,\nabla\sigma)=0$. Thus, applying (2.11) 
we conclude that
$$
	\varepsilon\bar\tau e^{2\sigma}-\tau = 0 \ .   \leqno (2.12)
$$
Then, (2.11) reduces to
$$
	Q=\nu \,d\sigma\otimes d\sigma    \leqno(2.11')
$$
or, according to (1.3) -
$$
	g(\nabla_X\nabla\sigma,Y)=(1-\nu)X\sigma Y\sigma \ .
$$
Hence, using the equation of Codazzi for $g$ and (2.6), we derive
$$
	(X\mu Y\sigma - Y\mu X\sigma)\,Z\sigma +
		\frac1n\{ X\,{\rm tr}\,h\,g(Y,Z)- Y\,{\rm tr}\,h\,g(X,Z) \}=0 \ .
$$
Here we assume that $Z$ is orthogonal to $\nabla\sigma$ and $X$, and $Y$ is 
not orthogonal to $Z$. The result is $X{\rm tr}\,h=0$. i.e. ${\rm tr}\,h$ is
a constant. Thus, by (2.8) and (2.10), $\tau$ and $\bar\tau$ are also constants.
If $\bar\tau \ne 0$, (2.12) implies $d\sigma=0$, which is a contradiction.
Let $\bar\tau =0$. According to (2.12), (2.8) and (2.10), $\tau={\rm tr}\,h=\lambda=0$.
By Lemma 3 c)
$$
	\overline R=e^{4\sigma}R  \leqno (2.13)
$$ 
On the other hand, from ${\rm tr}\,h=\lambda=0$ and (1.2), $(2.11')$, we obtain
$$
	\overline R=\varepsilon e^{2\sigma}R  \leqno (2.14)
$$
From (2.13) and (2.14) we find $(e^{2\sigma}-\varepsilon)\,R=0$ in $U_1$ and hence
this holds on $U$. Since $\sigma$ can not vanish in an open subset of $U$, it follows
$R=0$ in $U$, which proves our assertion.

\vspace{0.5in}

\centerline{\bf 3. PROOF OF THEOREM 1}

\vspace{0.3in}
First we assume that there exists a point $p$ of $M$ such that 
$\left\Arrowvert\nabla\sigma\right\Arrowvert^2\ne0$ in $p$. Then 
$\left\Arrowvert\nabla\sigma\right\Arrowvert^2\ne0$ 
in a neighbourhood $V$ of $p$. In (2.1) we assume that $X=Z=\nabla\sigma$ and that $Y$ 
is orthogonal to $\nabla\sigma$. Using (2.2), we obtain $Y\lambda=0$ in $V$. Hence 
$\nabla\lambda=\rho\nabla\sigma$ on $V$, where $\rho $ is a smooth function. Using again 
(2.1) with $X=\nabla\sigma$ and applying (2.2), we find
$$
	B=\rho \left\{ \frac1{\left\Arrowvert\nabla\sigma\right\Arrowvert^2}\,d\sigma\otimes d\sigma -\frac1n\,g \right\}
$$
in $V$. However, this contradicts the assumption that the set of nonquasi-umbilic points is dense.

So $\left\Arrowvert\nabla\sigma\right\Arrowvert^2=0$. Now, let us assume that $\nabla\sigma$ does not vanish 
at a point $p$ and hence, in an open set $U$. By Lemma 4 a) each point of $U$ is quasi-umbilic,
which is impossible.

Consequently $\nabla\sigma$ vanishes indentically in $M$, i.e. $\sigma$ is a constant. Then
$\lambda$ is also a constant. Indeed, assuming in (2.1) that $V$ is orthogonal to $X$ and that 
$Y=Z$, $g(Y,Y)\ne0$, we obtain $X\lambda=0$. 

Since $\sigma$ is a constant, (1.2) implies
$$
	f^*\overline R=\varepsilon e^{2\sigma}R \ .    \leqno (3.1)
$$
Let us assume that $f$ is not an isometry, i.e. $(\sigma,\varepsilon)\ne(0,1)$.
Then (3.1) and Lemma 3 c) yield
$$
	R=\alpha\varphi(h)+\beta\pi_1 \ ,   \leqno (3.2)
$$
where
$$
	\alpha=\frac{\lambda}{\varepsilon e^{-2\sigma}-1} \,, \qquad
	\beta=\frac{\lambda^2}{\varepsilon e^{-2\sigma}-1}
$$
are constants. From (3.2), by a standard way (see e,g, [5] or [6], Example 4), we
conclude that the Weyl conformal curvature tensor of $M$ vanishes identically. So,
if $n>3$ then $M$ is conformally flat, which is a contradiction. Let $n=3$. Using
(3.2) we find
$$
	S-\frac{\tau}4\,g=\alpha\,h+\frac{\beta}2\,g \ .    \leqno (3.3)
$$
Since $\alpha,\,\beta$ are constants, the equation of Codazzi and (3.3) imply (1.1).
Thus $M$ is conformally flat, which is not the case. Consequently $f$ is an isometry, 
i.e. $\sigma=0$, $\varepsilon=1$. Then by (3.1) and Lemma 3 c) we obtain
$$
	\lambda\varphi(h)+\lambda^2\pi_1=0 \ ,
$$
which implies
$$
	\lambda\{(n-2)h+g\,{\rm tr}\,h\} +(n-1)\lambda^2g=0 \ .   \leqno (3.4)
$$
But $M$ can not be totally umbilic. So (3.4) yields $\lambda=0$. Hence
$f_*\bar h=h$. Since $f$ is an isometry, this proves the theorem.

\vspace{0.7in}

\centerline{\bf 4. PROOF OF THEOREM 2}

\vspace{0.4in}
By Lemma 3 a), b) and (1.4) we conclude that $\lambda=0$. Putting
$X=\nabla\sigma$ in (2.1) and using (2.2), we obtain
$$
	\left\Arrowvert\nabla\sigma\right\Arrowvert^2B(Y,Z)=0 \ .
$$
Since the nonumbilic points are dense, this implies $\left\Arrowvert\nabla\sigma\right\Arrowvert^2=0$.
According to Lemma 4 b), $R=0$ in the open set $U$, in which $\nabla\sigma\ne0$. Thus
the point $p$, in which $R\ne0$, can not lie in the closure $\overline U$ of $U$. 
Consequently, the open set $M\setminus\overline U$ is nonempty. Note that $d\sigma=0$ in
$M\setminus\overline U$. Let $V$ be the connected component of $p$ in $M\setminus\overline U$.
Since $\sigma$ is a constant in $V$, (1.2) reduces to
$$
	f^*\overline R=\varepsilon e^{2\sigma}R     \leqno (4.1)
$$
in $V$. On the other hand, applying Lemma 3 c) with $\lambda=0$, we obtain
$$
	f^*\overline R=  e^{4\sigma}R \  .    \leqno (4.2)
$$
From (4.1) and (4.2) we find $(e^{4\sigma}-\varepsilon)R=0$. Since $p$ lies in $V$,
this implies $\sigma=0$  (in $V$) and $\varepsilon=1$. So we have $f^*\bar g=g$,
$f^*\bar h=h$ in $V$. Consequently, $f$ is a congruence of $V$ onto $f(V)$, which
completes the proof.

{\bf Remark.} If the manifolds in Theorem 2 are analytic or the set of points, in 
which $R$ is not zero, is dense, then $f$ is a congruence of $M$ onto $\overline M$.

\vspace{0.7in}
\centerline{\large REFERENCES}

\vspace{0.3in}

\noindent
\ 1. K\,u\,l\,k\,a\,r\,n\,i, R. S. Congruence of hypersurfaces. - Journ. Differ. Geom., {\bf 8},
1973, 

\ \ 95-102.

\noindent
\ 2. D\,a\,j\,c\,z\,e\,r, M., K. N\,o\,m\,i\,z\,u. On the boundedness of Ricci curvature of an indefinite 

\ \ metric. - Bol. Soc. Brasil. Mat., {\bf 11}, 1980, 25-30.
 
\noindent
\ 3. E\,i\,s\,e\,n\,h\,a\,r\,t, L. P.  Riemannian geometry. Princeton. University Press, 1949.

\noindent 
\ 4. K\,a\,s\,s\,a\,b\,o\,v, O. Diffeomorphisms of pseudo-Riemannian manifolds and the values of 

\ \ the curvature tensor on degenerate planes. Serdica, {\bf 15}, 1989, 78-86.

\noindent 
\ 5. C\,h\,e\,n, B.-Y., K. Y\,a\,n\,o. Manifolds with vanishing Weil or Bochner curvature tensor. 

\ \ - Journ. Math. Soc. Japan, {\bf 27}, 1975, 106-112.

\noindent
\ 6. N\,o\,m\,i\,z\,u, K. On the decomposition of the generalized curvature tensor. - 
Differential 

\ \ Geometry in honour of K. Yano. Kynokumya, Tokyo, 1972, 335-345.

\vspace {0.3in}

\ \ \ \ \ \ \ \ \ \ \ \ \ \ \
Received 30.03.1988

\end{document}